\documentclass[12pt]{amsart}
\usepackage{amssymb,amscd, fullpage}

\usepackage{palatino}
\usepackage{euler}

\numberwithin{equation}{section}

\newtheorem{theorem}{Theorem}[section]

\newtheorem{lemma}[theorem]{Lemma}
\theoremstyle{definition}

\theoremstyle{remark}
\newtheorem{remark}[theorem]{Remark}
\newtheorem{example}[theorem]{Example}

\newcommand{\cO}{\mathcal{O}}
\renewcommand{\L}{\mathcal{L}}
\newcommand {\cL}{\mathcal{L}}

\newcommand{\R}{\mathbb{R}}
\newcommand{\C}{\mathbb{C}}
\newcommand{\Z}{\mathbb{Z}}

\newcommand{\T}{\mathbb{T}}
\newcommand\lie[1]{\mathfrak{#1}}

\newcommand{\fh}{\lie{h}}
\newcommand{\fg}{\lie{g}}

\newcommand{\fk}{\lie{k}}

\def    \inv    {^{-1}}

\newcommand\cone{\mathaccent23}

\newcommand\Span{\mathop{\rm span}\nolimits}
\newcommand\affspan{\mathop{\rm affspan}\nolimits}

\newcommand{\intF} {{\mathaccent23{F}}}

\begin{document}

\title{K\"ahler metrics on singular toric varieties }

\author{D. Burns \and V. Guillemin \and E. Lerman}\thanks{Supported in
part by NSF grants DMS-0104047 and DMS-0514070 (DB), DMS-0104116 and DMS-0408993 (VG), and
DMS-0204448(EL)} 
\address{University of Michigan, Ann
Arbor, MI 48109}

\email{dburns@umich.edu}

\address{M.I.T., Cambridge, MA 02139}

\email{vwg@math.mit.edu}

\address{University of Illinois, Urbana, IL 61801 } 
\email{lerman@math.uiuc.edu}


\begin{abstract}
We extend Guillemin's formula for K\"ahler potentials on toric
manifolds to singular quotients of $\C^N$ and $\C P^N$.
\end{abstract}
\maketitle

\section{Introduction}

Let $G$ be a torus with Lie algebra $\fg$ and integral lattice $\Z_G
\subset \fg$.  Let $u_1, \ldots u_N \in \Z_G$ be a set of primitive vectors 
which span $\fg$ over $\R$.  Let $\lambda _1, \ldots,
\lambda _N \in \R$ and let
\[
P = P_{u,\lambda} := 
\{ \eta \in\fg^*\mid \langle\eta , u_j\rangle -\lambda_j \geq 0,
\quad 1\leq j \leq N\}
\]
be the corresponding polyhedral set.  We \emph{ assume that $P$ has a
non-empty interior and that the collection of inequalities defining
$P$ is minimal}: if we drop the condition that $ \langle\eta ,
u_j\rangle -\lambda_j \geq 0$ for some index $j$ then the resulting
set is strictly bigger than $P$.  

A well-known construction of Delzant, suitably tweaked, produces a
symplectic stratified space $M_P$   with an effective Hamiltonian
action of the torus $G$ and associated moment map $\phi =\phi_P: M_P\to
\fg^*$ such that $\phi  (M_P) = P$.  We will review the construction 
below.  The space $M_P$ is a symplectic quotient of $\C^N$ by a 
compact abelian subgroup $K$ of the standard torus $\T^N$.  Therefore,
by a theorem of Heinzner and Loose $M_P$ is a complex analytic space
\cite{HL}.  Moreover $M_P$ is a K\"ahler space \cite[(3.5)]{HL} and
\cite{HHuL}.  Even though in general the space $M_P$ is singular, the
preimages of open faces of $P$ under the moment map $\phi_P$ are
smooth K\"ahler manifolds.  The main result of the paper are formulas
for the K\"ahler forms on these manifolds.
In particular 
we will show that the K\"ahler form $\omega$ on the preimage
$\phi_P\inv(\mathaccent23P)$ of the interior  $\mathaccent23P$ of
the polyhedral set $P$ is given by the formula~(\ref{theformula}) below:
\begin{equation}
\label{theformula}
  \omega = \sqrt{-1} \partial \bar{\partial} \phi^* (
\sum_{j=1}^N \lambda _j \log ( u_j -\lambda _j) + u_j ),
\end{equation}
where we think of $u_j \in \Z_G$ as a function on $\fg^*$.

Formula (\ref{theformula}) was originally proved by Guillemin in the
case where $M_P$ is a compact manifold (and thus $P$ is a simple
unimodular polytope, also known as a Delzant polytope).  It was
extended to the case of compact orbifolds by Abreu \cite{Abr}.
Calderbank, David and Gauduchon gave two new proofs of Guillemin's
formula (for orbifolds) in \cite{CDG}.  One of their proofs was
simplified further in
\cite{BG}.

As we just mentioned, for generic values of $\lambda$ the polyhedral
set $P$ is simple and consequently $M_P$ is at worse an orbifold.  But
for arbitrary values of $\lambda$ it may have more serious
singularities.  Of particular interest is the singular case where $P$
is a cone on a simple polytope.  Then there is only one singular
point, and the link of the singularity is a Sasakian orbifold.  Such
orbifolds, especially the ones with Sasaki-Einstein metrics, have
attracted some attention in string theory.  They play a role in the
AdS/CFT correspondence \cite{MS}.

If the polyhedral set $P$ is a polytope, i.e., if $P$ is compact, then
as a symplectic space $M_P$ may also be obtained as a symplectic
quotient of $\C P^N$.  In this case the Fubini-Study form on $\C P^N$
will induce a K\"ahler structure on $M_P$, which is {\em different}
from the one induced by the flat metric on $\C^N$ even in the case
where $M_P$ is smooth.  We will give a formula for this K\"ahler
structure as well.

\subsection*{Acknowledgments}  We thank Peter Heinzner for patiently
 answering our queries.

\section{The ``Delzant'' construction: toric varieties as K\"ahler quotients}

It will be convenient for us to fix the following notation.  As in the
introduction, let $G$ be a torus with Lie algebra $\fg$ and integral
lattice $\Z_G
\subset \fg$.  Let $u_1, \ldots u_N \in \Z_G$ be a set of primitive vectors 
which span $\fg$ over $\R$.  Let $\lambda _1, \ldots,
\lambda _N \in \R$ and let
\begin{equation} 
\label{eqP}
P = P_{u,\lambda} := 
\{ \eta \in\fg^*\mid \langle\eta , u_j\rangle -\lambda_j \geq 0,
 \quad 1\leq j\leq N \}
\end{equation}
be the corresponding polyhedral set.  As above {\em we assume that
$P$ has the non-empty interior and that the collection of inequalities
defining $P$ is minimal}.   Let $A:\Z^N \to \Z_G$ be the $\Z$-linear map 
given by
\[
A (x_1, \ldots, x_N) = \sum x_i u_i .
\]
That is, $A$ is defined by sending the standard basis vector $e_i$ of
$\Z^N$ to $u_i$.  Let $A$ also denote the $\R$-linear extension $\R^N
\to \fg$.  Let $\fk = \ker A$ and let $B:\fk \to \R^N$ denote the
inclusion.  The map $A$ induces a surjective map of Lie groups 
\[
\bar {A} : \T ^N = \R^N/\Z^N \to \fg/\Z_G = G. 
\]
Let $K = \ker \bar {A}$ and let $\bar{B} :K \to \T^N$ denote the
corresponding inclusion.  The group $K$ is a compact abelian group
which need not be connected.  It's easy to see that the Lie algebra of
$K$ is $\fk$.

We have a short exact sequence of abelian Lie algebras:
\[
0 \to \fk \stackrel{B}{\to } \R^N \stackrel{A}{\to } \fg \to 0.
\]
Let 
\[
0 \to \fg^* \stackrel{A^*}{\to } (\R^N)^* \stackrel{B^*}{\to } \fk^* \to 0
\]
be the dual sequence.  Note that $\ker B^* = A^* (\fg^*) = \fk^\circ$
where $\fk^\circ$ denotes the annihilator of $\fk$ in $(\R^N)^*$.  Let
$\{e_i^*\}$ denote the dual basis of $(\R^N)^*$ and let $\lambda = \sum
\lambda _i e_i^*$.  We note that 
\[
(B^*)\inv (B^* (- \lambda)) = -\lambda + \fk^\circ = -\lambda + A^* (\fg^*).
\]
In particular $(B^*)\inv (B^* (- \lambda))$ is the image of the affine
embedding
\begin{equation}
\label{eq.affine}
\iota_\lambda : \fg^* \hookrightarrow (\R^N)^*, \quad
\iota_\lambda (\ell) = -\lambda + A^* (\ell).
\end{equation}

\begin{lemma} 
\label{lem21}
 Let $P$ be the polyhedral set defined by (\ref{eqP}) above.
\begin{enumerate}
\item There exists a K\"ahler space $M_P$ with an effective
holomorphic Hamiltonian action of the torus $G$ so that the image of
the associated moment map $\phi_P:M_P \to \fg^*$ is $P$.
\item For every open face $\cone F$, the preimage $\phi_P\inv (\cone
F)$ is the K\"ahler quotient of a complex torus $(\C^{\times})^{N_F}$
by a compact subgroup $K_F$ of the compact torus $\T^{N_F} \subset
(\C^{\times})^{N_F}$.  Here the number $N_F$ and the group $K_F$
depend on the face $F$.

\item If the set $P$ is bounded, then $M_P$ can also be constructed as
a K\"ahler quotient of $\C P^N$.
\end{enumerate}
\end{lemma}

\begin{proof}
For every index $i$ and any $\eta \in \fg^*$
\begin{multline}
\langle \eta, Ae_i\rangle - \lambda_i = \langle A^* \eta, e_i\rangle -
\langle \sum \lambda_j e_j^*, e_i\rangle =
\langle A^* \eta - \lambda, e_i\rangle =
 \langle \iota_\lambda (\eta), e_i\rangle.
\end{multline}
Therefore 
\[
\iota_\lambda (P) = \{ \ell \in (\R^N)^* \mid 
\langle \ell, e_i \rangle \geq 0, \, 1\leq i \leq N\} \cap 
\iota_\lambda (\fg^*).
\]
More generally, if $\mathaccent23F \subset P$ is an open face, there
is a unique subset $I_F = I\subset \{1, \ldots, N\}$ so that
\begin{equation}
\label{eq2.4}
\mathaccent23F = \bigcap _{j\not \in I}
 \{ \eta \in \fg^* \mid \langle \eta, u_j\rangle - \lambda _j > 0\} 
\cap\bigcap _{j\in I}
\{ \eta \in \fg^* \mid \langle \eta, u_j\rangle - \lambda _j = 0\} .
\end{equation}
Therefore
\begin{equation}
\label{eq2.5}
\iota_\lambda (\mathaccent23F) = 
\iota_\lambda (\fg^*) \cap \bigcap_{j\not \in I}
 \{ \ell \in (\R^N)^* \mid \langle \ell, e_j\rangle > 0\} 
\cap \bigcap_{j\in I}
 \{ \ell \in (\R^N)^* \mid \langle \ell, e_j\rangle = 0\} .
\end{equation}
The moment map $\phi$ for the action of $\T^N$ on 
$(\C^N, \sqrt{-1}\sum dz_j \wedge d\bar{z}_j)$  is given by
\[
\phi (z) = \sum |z_j|^2 e_j^* .
\]
Hence 
\[
\phi (\C^N) = \{ \ell \in (\R^N)^* \mid 
\langle \ell, e_i \rangle \geq 0, \quad
1\leq i \leq N\}.
\]
The moment map $\phi_K$ for the action of $K$ on $\C^N$ is the composition 
\[
\phi_K = B^* \circ \phi.
\]
Let $\nu = B^* (-\lambda)$.  We argue that 
\[
\phi (\phi_K \inv (\nu)) = \iota_\lambda (P). 
\]
Indeed,
\begin{equation}
\begin{split}
\phi_K \inv (\nu) &= \phi\inv ( (B^*)\inv (\nu)) \\
	&= \phi\inv ( (B^*)\inv (B^* (-\lambda)))\\
	&=\phi \inv (\iota_\lambda (\fg^*))\\
	&= \phi\inv (\phi (\C^N) \cap \iota _\lambda (\fg^*))\\
	&= \phi \inv (\iota _\lambda (P)) .
\end{split}
\end{equation}
Therefore
\[
\phi (\phi_K \inv (\nu)) = \iota _\lambda (P).
\]
The restriction $\phi|_{\phi_K \inv (\nu)}$ descends to a map
$\bar{\phi} : M_P \equiv \phi_K \inv (\nu)/K \to \iota_\lambda
(\fg^*)$.  It is not hard to see that the composition $\phi_P$ of
$\bar{\phi}$ with the isomorphism $\iota_\lambda (\fg^*)
\stackrel{\simeq}{\to } \fg^*$ is a moment map for the action of $G$
on the symplectic quotient (symplectic stratified space) $M_P$.  Since
the isomorphism $\iota_\lambda (\fg^*) \to \fg^*$ obviously maps
$\iota_\lambda (P)$ to $P$, we conclude that the image of $\phi_P: M_P
\to \fg^*$ is exactly $P$. This proves (1).\\

To prove (2) we  define a bit more notation.
For a subset $I\subset \{1,\ldots, N\}$ we define the corresponding
coordinate subspace
\[
V_I := \{ z\in \C^N \mid j\in I \Rightarrow z_j = 0\}.
\]
Its ``interior'' $\mathaccent23V _I$ is defined by
\[
\mathaccent23V _I:= \{ z\in \C^N \mid j\in I \Leftrightarrow z_j = 0\}.
\]
Also, let 
\[
\T^N_I := \{ a\in \T^N \mid j\not \in I \Rightarrow a_j =1\} .
\]
The sets $V_I$, $\mathaccent23V _I$ are K\"ahler submanifolds of
$\C^N$ preserved by the action of $\T^N$.  They are both fixed by
$\T^N_I$ with $\mathaccent23V _I$ being precisely the set of points of
orbit type $\T^N_I$.

The restriction $\phi_K |_{\mathaccent23V _I}$ is a moment map for the
action of $K$ on $\mathaccent23V _I$.  Moreover, for any $\nu\in
\fk^*$
\[
\phi_K\inv (\nu) \cap \mathaccent23V _I
= (\phi_K|_{\mathaccent23V _I})\inv (\nu).
\]
Hence
\[
(\phi_K\inv (\nu) \cap \mathaccent23V _I)/K
= (\phi_K|_{\mathaccent23V _I})\inv (\nu)/K.
\]
While the action of $K$ on $\mathaccent23V _I$ need not be free, the
action of
\[
K_I := K/(K\cap \T^N_I)
\]
on $\mathaccent23V _I$ {\em is} free.  Therefore, the quotient
$(\phi_K\inv (\nu) \cap \mathaccent23V _I)/K$ may be interpreted as a
{\em regular} K\"ahler quotient of $\mathaccent23V _I$ by the
Hamiltonian action of $K_I$:
\begin{equation}
\label{eq2.7}
(\phi_K\inv (\nu) \cap \mathaccent23V _I)/K = 
\mathaccent23V _I/\!/_{\nu_I} K_I
\end{equation}
for an appropriate value $\nu_I \in \fk_I^*$ of the $K_I$ moment map.

Given a face $F$, let $I= I_F$ be the corresponding subset of $\{1,
\ldots, N\}$.  Then, by (\ref{eq2.5}),
\begin{multline}
\{z\in \C^N \mid \phi (z) \in \iota_\lambda ( \mathaccent23F)\} =\\
\{z\in \C^N \mid \phi (z) \in \iota_\lambda (\fg^*), \langle \phi (z),
e_j\rangle > 0 \text{ for }j\not\in I, \langle \phi (z), e_j\rangle =0
\text{ for }j\in I\}\\
= \phi_K \inv (\nu) \cap  \mathaccent23V _I.
\end{multline}
Therefore,
\[
\phi_K \inv (\nu) \cap  \mathaccent23V _I
 = \phi\inv (\iota_\lambda ( \mathaccent23F)).
\]
It follows from the definition of $\phi_P$ that
\[
(\phi_K \inv (\nu) \cap  \mathaccent23V _I)/K
 = \phi_P \inv  ( \mathaccent23F).
\]
By (\ref{eq2.7}) we conclude that 
\begin{equation}
\label{eq-*}
\phi_P \inv  ( \mathaccent23F) = \mathaccent23V _I/\!/_{\nu_I} K_I .
\end{equation}
This proves (2).\\

If $P$ is compact, then $\iota_\lambda (P) \subset (\R^N)^*$ is
bounded.  Hence $\iota_\lambda (P)$ is contained in a sufficiently
large multiple of the standard simplex.  Any such simplex is the image
of $\C P^N$ under the moment map for the standard action of $\T^N$
with the K\"ahler form on $\C P^N$ being the appropriate multiple of
the standard Fubini-Study form.  This proves (3).
\end{proof}

\begin{remark}
It follows from the results of Heinzner and his collaborators
\cite{H}, in particular from \cite{HHu}, that the action of $G$ on
$M_P$ extends to an action the complexified group $G^\C$.  This action
of $G^\C$ has a dense open orbit.  In other words $M_P$ is a {\em  toric}
K\"ahler space.
\end{remark}

\section{K\"ahler potentials, Legendre transforms and symplectic
quotients}

We start this section by recalling a result of Guillemin (Theorem~4.2
and Theorem~4.3 in \cite{G}):
\begin{lemma}
\label{lem2.3}
Suppose the action of $\T^N$ on $(\C^{\times})^N = \R^N \times \sqrt{-1} \T^N$
preserves a K\"ahler form $\omega$ and is Hamiltonian.  Then there
exists a $\T^N$-invariant function $f$ on $(\C^{\times})^N$ such that
$\omega = i\partial \bar{\partial} f$.  Additionally
\[
{\mathcal L}_f \circ \pi : (\C^{\times})^N \to (\R^N)^*
\]
is a moment map for the action of $\T^N$ on $((\C^{\times})^N,
\omega)$.  Here $\pi: \R^N \times \sqrt{-1}\T^N \to \R^N$ is the projection
and ${\mathcal L}_f: \R^ N \to (\R^N)^*$ is the Legendre transform of
$f$, where we have identified $f\in C^\infty ((\C^{\times})^N)^{\T^N}$
with a function on $\R^N$.

The same result holds with $(\C^{\times})^N$ replaced by $U\times
\sqrt{-1}\T^N$ for any contractible open set $U\subset \R^N$.
\end{lemma}

\begin{lemma}
\label{lem2.6}
Let $f:V \to \R$ be a (strictly) convex function on a finite
dimensional vector space $V$, let $A:W \to V$ be an injective linear
map, $x\in V$ be a point and 
\[
j: W \to V, \quad j(w) = Aw +x
\]
an affine map.  Then $f\circ j: W \to \R$ is (strictly) convex and the
associated Legendre transform ${\mathcal L}_{f\circ j} :W \to W^*$ is
given by
\[
{\mathcal L}_{f\circ j}  = A^* \circ {\mathcal L}_f \circ j,
\]
where $A^*: V^* \to W^*$ is the dual map.
\end{lemma}
\begin{proof}
By the chain rule and the definition of the Legendre transform,
${\mathcal L}_{f\circ j} (w) = d (f\circ j)_w = df_{j(w)} \circ dj_w =
{\mathcal L}_f (j(w)) \circ A = A^* \circ {\mathcal L}_f \circ j (w)$
for any $w\in W$.  
\end{proof}

\begin{lemma} 
\label{lem-potential}
Let $f\in C^\infty (\R^N)$ be a strictly convex function and $\omega =
\sqrt {-1} \partial \bar{\partial}\, \pi_N^* f$ the corresponding
$\T^N$-invariant K\"ahler form on $(\C^{\times})^N = \R^N
\times\sqrt{-1} \T^N$
(here $\pi_N:(\C^{\times})^N \to \R^N $ is the projection). 
Let $\phi = {\mathcal L}_f \circ \pi_N : (\C^{\times})^N \to (\R^N)^*$
denote the associated moment map.

Let $K\subset \T^N$ be a closed subgroup and let $G= \T^N/K$.  For any
$\nu\in \fk^*$ the symplectic quotient $(\C^{\times})^N/\!/_\nu K$ is
biholomorphic to $U\times \sqrt{-1}G \subset \fg \times \sqrt{-1} G = G^\C$ 
where $U\subset \fg$ is an open contractible set.  Hence the reduced K\"ahler
 form $\omega_\nu$ has a potential $f_\nu$.

Moreover, the Legendre-Fenchel dual $f_\nu^*$ of the K\"ahler potential 
$f_\nu$  is given by
\begin{equation} 
\label{eq-reduced-pot}
f_\nu^* = f^* \circ \iota_\lambda
\end{equation}
where 
$\iota_\lambda : \fg^* \to (\R^N)^*$ is the affine embedding
(\ref{eq.affine}) and $-\lambda $ is a point in $(B^* ) \inv (\nu)$.
\end{lemma}

\begin{proof}
It is no loss of generality to assume that the group $K$ is
connected. Then $\T^N \simeq K \times G$.  Consequently $\R^N\simeq
\fk \times \fg$ and the short exact sequence
\[
0 \to \fk \stackrel{B}{\to } \R^N \stackrel{A}{\to } \fg \to 0.
\]
splits.  Let 
\[
\pi_K : \R^N \to \fk  \quad \text{and} \quad \iota_\fg : \fg \to \R^N
\]
denote the maps defined by the splitting.  The moment map $\phi_K :
(\C^{\times})^N \to \fk^*$ for the action of $K$ on $((\C^{\times})^N,
\omega)$ is the composition
\[
\phi_K = B^* \circ \phi = B^* \circ {\mathcal L}_f \circ \pi_N .
\]

Let 
\[
\Delta = (B^*)\inv (\nu) \cap \phi ((\C^{\times})^N) =
(B^*)\inv (\nu) \cap 
{ \mathcal L}_f (\R^N).
\]
Then $\Delta$ is the intersection of an affine hyperplane with a
convex set, hence is contractible.

Since the action of $K$ on $\phi_K \inv (\nu)$ is free, $K^\C \cdot
\phi_K\inv (\nu)$ is an open subset of $(\C^{\times})^N$ and $K^\C$
acts freely on it.  Moreover, for each $x\in \phi_K\inv (\nu)$ the
orbit $K^\C \cdot x$ intersects the level set $\phi_K\inv (\nu)$
transversely and
\[
K^\C \cdot x \cap \phi_K\inv (\nu) = K\cdot x
\]
(See \cite[pp.\ 526--527]{GS}).
It follows that the restriction
\[
\pi_K|_{{\mathcal L}_f\inv (\Delta)}: {{\mathcal L}_f\inv (\Delta)}
\to \fk
\]
is 1-1 and a local diffeomorphism.   Hence 
\[
U = \pi_K ({\mathcal L}_f\inv (\Delta))
\]
is a contractible open set.

On the other hand, the restriction $\omega|_{\phi_K\inv (\nu)}$
descends to a K\"ahler form $\omega_\nu$ on the symplectic quotient
\[
(\C^{\times})^N) /\!/_\nu K := \phi_K\inv (\nu)/K.
\]
Moreover, since $\omega$ is $\T^N$ invariant, $\omega_\nu$ is
$G$-invariant.  Note that
\[
(\C^{\times})^N) /\!/_\nu K \simeq U\times \sqrt{-1} G \subset G^\C.
\]
By Lemma~\ref{lem2.3} there exists $f_\nu\in C^\infty (U)$ such that
\[
\omega_\nu = \sqrt{-1} \partial \bar{\partial} f _\nu.
\]
The potential $f_\nu$ defines a moment map 
\[
\phi_G : U\times \sqrt{-1} G \to \fg^* 
\]
with
\[
  \phi_G = {\mathcal L}_{f_\nu} \circ \pi_G
\]
where $\pi_G : U\times \sqrt{-1} G \to U$ is the projection.
Moreover, by adjusting $f_\nu$ \cite{BG} we may arrange for the diagram
\begin{equation}
\label{cd1}
\begin{CD}
\phi_K\inv (\nu) @>\phi >> \Delta \subset (\R^N)^*\\
@V /KVV	   @AA \iota_\lambda A\\
U \times \sqrt{-1} G  @>> \phi_G > \fg^*
\end{CD}
\end{equation}
to commute.  That is, the moment map $\phi_G$ is defined up to a
constant and the potential $f_\nu$ is defined up to a pluri-harmonic
$G$-invariant function.  By adding an appropriate pluri-harmonic
function to $f_\nu$ we can change $\phi_G$ by any constant we want.
Since $\phi = {\mathcal L}_f \circ \pi_N$ and since $\phi_G =
\L_{f_\nu} \circ \pi_G$, it follows from (\ref{cd1}) that the diagram
below commutes as well:
\begin{equation}
\label{cd2}
\begin{CD}
{\L}_f \inv (\Delta) @> \L_f >> \Delta \subset (\R^N)^*\\
@V \pi_K VV	   @AA \iota_\lambda A \\ 
U   @>> \L_{f_\nu} > \fg^*
\end{CD} .
\end{equation}
Since $(\L_{f_\nu})\inv = \L_{f^*_\nu}$, where $f^*_\nu$ is the
Legendre-Fenchel dual of $f_\nu$,
\[
\L_{f^*_\nu} = \pi_K \circ (\L _f)\inv \circ \iota_\lambda =
\pi_K \circ (\L _{f^*})\circ \iota_\lambda.
\]
By Lemma~\ref{lem2.6},
\[
\L_{f^*_\nu} = \L_{f^*} \circ \iota_\lambda .
\]
Therefore, up to a constant,
\[
f^*_\nu = f^* \circ \iota_\lambda.
\]
\end{proof}

\section{From potentials to  dual potentials and back again}

We start by making two observations.  Let $V$ be a real finite
dimensional vector space, $V^*$ its dual, $\cO \subset V$ an open
set, $\varphi\in C^\infty (\cO)$ a strictly convex function,
$\L_\varphi:\cO\to V^*$ the Legendre transform (which we assume to be
invertible), $\cO^* = \L_\varphi (\cO)$ and $\varphi^* \in C^\infty
(\cO^*)$ the Fenchel dual of $\varphi$.

\begin{lemma}
\label{lem3.1} Under the above assumptions,
$\varphi = ({\L_\varphi})^*h$ where $h:\cO^*\to \R$ is given by
\[
h(\eta) = \langle \eta, (d\varphi^*)_\eta \rangle - \varphi^* (\eta)
\]
where we think of $(d\varphi^*)_\eta \in T^*_\eta \cO^* $ as an
element of $(V^*)^* = V$.
\end{lemma}
\begin{proof}
By definition of the Fenchel dual
\[
\varphi (s) + \varphi^* (\eta ) = \langle \eta, s\rangle
\]
for $\eta = \L _{\varphi} (s)$.
Hence
\[
\varphi (s) = \langle \eta, s\rangle -  \varphi^* (\eta ) =\langle
\eta, (\cL_\varphi)\inv (\eta) \rangle -  \varphi^* (\eta ) =
\langle \eta, \cL_{\varphi^*} (\eta) \rangle -  \varphi^* (\eta )
\]
and the result follows since $\cL_{\varphi^*} (\eta) = (d\varphi^*)_\eta$.
\end{proof}

\begin{lemma}
\label{lem3.2}
We keep the above notation.  Suppose additionally that the dual
potential $\varphi^*$ has the following special form:
\[
\varphi^* (\eta)  = \sum_{i=1}^N f_i (u_i (\eta)-\lambda_i),
\]
where $u_1, \ldots, u_N$ are vectors in $V$ (thought of as linear
functionals $u_i:V^*\to \R$), $\lambda_i\in \R$ are constants and
$f_i$'s are functions of one variable.  Then
\begin{equation}\label{eq3.1}
h (\eta) = \sum _{i=1}^N 
\Big ( f_i' (u_i (\eta)-\lambda_i)\;u_i (\eta)
 -f_i (u_i (\eta)-\lambda_i)\Big ).
\end{equation}
\end{lemma}

\begin{proof} Observe that
$d(f_i \circ (u_i-\lambda_i)) _\eta = f_i' (u_i (\eta)-\lambda_i)
d(u_i-\lambda_i)_\eta = f_i' (u_i (\eta)-\lambda_i)\; u_i$ since $u_i$
is linear.  Hence
\[
\langle \eta, (d\varphi^*)_\eta\rangle = \langle \eta, \sum f_i' (u_i
(\eta)-\lambda_i) u_i\rangle = 
 \sum f_i' (u_i (\eta)-\lambda_i)\; u_i (\eta)
\]
and (\ref{eq3.1}) follows from Lemma~\ref{lem3.1} above.
\end{proof}

\begin{example}
\label{ex3.6}
We use the lemma above to argue that for the standard action of $\T^N$
on $(\C^N, \sqrt{-1}\partial \overline{\partial} ||z||^2)$, the dual
potential $\varphi^*$ is given by
\[
\varphi^* = \sum_{i=1}^N e_i \log e_i,
\]
where $e_1, \ldots, e_N$ is the standard basis of $\R^N = Lie (\T^N)$.

Indeed, the homogeneous moment map $\Phi: \C^N\to (\R^N)^*$ for the
standard action of $\T^N$ is given by
\[
\Phi (z) = \sum |z_j|^2 e_j^*,
\]
where $\{e_j^*\}$ is the basis dual to $\{e_j\}$.  Hence 
\[
 ||z||^2  = \Phi^* (\sum e_j).
\]
On the other hand, if $\varphi^* = \sum e_j \log e_j$, then
\[
 \varphi^* = \sum f\circ e_j
\]
where $f(x) = x\log x$.  Since $f'(x) = \log x + 1$, equation
(\ref{eq3.1}) becomes
\[
 h = \sum (\log e_j + 1)e_j -\sum e_j \log e_j  = \sum e_j.
\]
Therefore, $\varphi^* = \sum e_j\log e_j$ is, indeed, the dual potential.
\end{example}

We are now in position to prove (\ref{theformula}).  
\begin{theorem}
Let $G$ be a torus, $P\subset \fg^*$ the polyhedral set defined by
(\ref{eqP}), $M_P = \C^N/\!/_\nu K$ the K\"ahler $G$-space with moment
map $\phi_P:M_P \to \fg^*$ constructed in Lemma~\ref{lem21}~(1).  Then the
K\"ahler form $\omega_P$ on $\mathaccent23M_P: = \phi_P\inv
(\mathaccent23P)$ is given by:
\[
  \omega_P = \sqrt{-1} \partial \bar{\partial} \phi_P ^* (
\sum_{j=1}^N \lambda _j \log ( u_j -\lambda _j) + u_j ),
\]
\end{theorem}
\begin{proof}
By Lemma~\ref{lem21}, $\mathaccent23M_P = (\C^\times)^N/\!/_\nu K$
where $K\subset \T^N$ is a closed subgroup.  By
Lemma~\ref{lem-potential} the dual potential $\varphi_P^*$ on
$\mathaccent23P$ is given by
\[
\varphi_P^* = \varphi^* \circ \iota_\lambda
\]
where $\varphi^*$ is the potential on the open orthant in $(\R^N)^*$
dual to the flat metric potential $\varphi (z) = ||z||^2$ on
$(\C^\times)^N$.  By Example~\ref{ex3.6} $\varphi^* = \sum e_j \log
e_j$.  Since $\iota_\lambda^* e_j = u_j -\lambda_j$, 
\begin{equation} \label{new2.13}
\varphi _P^*  = \sum (u_j - \lambda_j) \log (u_j - \lambda_j).
\end{equation}
By Lemmas~\ref{lem3.1} and \ref{lem3.2}, the potential $\varphi_P$ is
given by 
\[
\varphi_P = \phi_P^* h
\]
where 
\[
h = \sum (\log (u_j -\lambda_j) +1) u_j
 - \sum (u_j -\lambda_j)\log (u_j -\lambda_j),
\]
(c.f. (\ref{eq3.1}).
Therefore
\[
\varphi_P = \phi_P^*\left( \sum_{i=1}^N (\lambda_j \log (u_j
-\lambda_j) + u_j)
\right)
\]
and we are done.
\end{proof}

\section{K\"ahler potentials on the preimages of faces}

Once again let $P\subset \fg^*$ be a polyhedral set given by
(\ref{eqP}). Recall that in section~2 we canonically associated to
this set a K\"ahler quotient $M_P$ of $\C^N$ which carries an
effective holomorphic and Hamiltonian action of the torus $G$ with a
moment map $\phi_P:M_P\to \fg^*$.
 Let $F\subset P$ be a face.  Its interior $\intF$ is given by:
\[
\mathaccent23F = \bigcap _{j\not \in I}
 \{ \eta \in \fg^* \mid \langle \eta, u_j\rangle - \lambda _j > 0\} 
\cap\bigcap _{j\in I}
\{ \eta \in \fg^* \mid \langle \eta, u_j\rangle - \lambda _j = 0\}
\]
for some nonempty subset $I$ of $\{1, \ldots, N\}$.  We have seen in
the proof of Lemma~\ref{lem21} that the preimage 
\[
M_\intF :=
\phi_P\inv (\intF)
\]
 is the K\"ahler quotient of $\cone V_I$ by a
compact abelian group $K_I$.  Therefore there is a potential
$\varphi_F^* \in C^\infty (\cone F)$ dual to the K\"ahler potential
$\varphi_F$ on $M_\intF$.  The goal of this section is to compute the
dual potential $\varphi_F^* $ ``explicitly.'' Lemmas~\ref{lem3.1} and
\ref{lem3.2} will then give us an analogue of (\ref{theformula}) for
the K\"ahler metric on $M_\intF$.

The K\"ahler potential $\varphi_I$ on
$\cone V_I$ for the flat metric induced from $\C^N$ is given by
\[
\varphi_I (z) = \sum _{j\not\in I} |z_j|^2.
\]
The restriction of the moment map $\phi:\C^N \to (\R^N)^*$ to $\cone
V_I$ is a moment map for the action of the torus 
\[
H_I := \T^N/\T^N_I.
\]
Note that
\[
\phi (\cone V_I) = \{\sum_{i\not \in I} a_i e_i^* \mid a_i>0\} .
\]
This set is an open subset in $\Span _{i\not \in I} \{ e_j^*\} \simeq
\fh_I^*$.   From now on we identify $\fh_I^*$ with 
$\Span _{i\not \in I} \{ e_j^*\}$.  The dual potential $\varphi_I^*\in
C^\infty (\phi (\cone V_I))$ is easily seen to be 
\[
\varphi_I^* = \sum_{j\not \in I} e_j\log e_j.
\]
The manifold $M_\intF$ is a Hamiltonian $G$ space, but the group $G$
doesn't act effectively.  So we cannot yet apply
Lemma~\ref{lem-potential} as we would like.  Let $G_I$ denote the
quotient of $G$ that does act effectively on $M_\intF$.  It is 
 isomorphic to the quotient
$H_I/K_I$. The dual of its Lie algebra $\fg_I^*$ is naturally embedded
in $\fg^*$:
\[
\fg_I^* = \{\eta\in \fg^* \mid \langle \eta, u_i\rangle = 0 \textrm{
for all } i\in I\}.
\]
Note also that the affine span $\affspan \intF$ of $\intF \subset
\fg^*$ is the translation of $\fg_I^*$ by an element $\eta_0\in
\intF$, as it should be.  Let $\gamma_I:\fg_I^* \to \affspan \intF
\subset \fg^*$ denote the affine embedding.  Then there exists an
affine embedding $\iota_I:\fg_I^* \hookrightarrow \fh_I^*$ so that the
diagram
\begin{equation}
\label{eq5.1}
\begin{CD}
\fh_I^* @> >>  (\R^N)^*\\
@A \iota_I AA	   @AA \iota_\lambda A\\
\fg_I^*  @>> \gamma_I > \fg^*
\end{CD}
\end{equation}
commutes.  Here the top arrow simply identifies $\fh_I^*$ with
$\Span_{i\not \in I}\{e_i^*\}$.  Since $\gamma_I$ is an embedding, we
may think of $\varphi_F^*$ as living on $\intF \subset \gamma_I
(\fg_I^*)$.
Therefore, by Lemma~\ref{lem-potential}, 
\begin{equation}
  \varphi_F^* = (\varphi_I^* \circ \iota_\lambda)|_\intF.
\end{equation}
Let 
\[
v_j = u_j |_\intF.
\]
These functions are affine, but not necessarily linear.
Then 
\[
  (e_j\circ \iota_\lambda)|_\intF = (u_j -\lambda_j)|_\intF = v_j - \lambda_j.
\]
Therefore
\[
  \varphi_F^* = (\varphi_I^* \circ \iota_\lambda)|_\intF =
\sum _{j\not \in I} (v_j -\lambda_j) \log (v_j -\lambda_j).
\]

To get a nicer formula for the potential on $M_\intF$ we now make a
simplifying assumption, namely, that $0\in \intF.$ Then $v_j = u_j|_{\fg_I^*}$ and, in particular, it is
{\em linear} for all $j$.  Hence Lemmas~\ref{lem3.1} and \ref{lem3.2}
apply, and we obtain:

\begin{theorem}
Under the simplifying assumption above, the K\"ahler form $\omega_F$
on $M_\intF$ is given by
\[
\omega_F = \sqrt{-1} \partial \bar{\partial} (\phi_P|_{M_\intF}) ^* (
\sum_{j\not \in I} \lambda _j \log ( v_j -\lambda _j) + v_j )
\]
\end{theorem}

Alternatively we may take the isomorphism $\gamma_I:\fg_I^* \to
\affspan \intF$ explicitly into account and think of $\varphi_F^*$ as
living on an open subset of $\fg_I^*$.  Then, by Lemma~\ref{lem-potential},
\[
 \varphi_F^* = \varphi_I^* \circ \iota_\lambda \circ \gamma_I. 
\]
Since
\[
 e_i \circ \iota_\lambda \circ \gamma_I = u_i|_{\fg^*_I} + u_i
(\eta_0) - \lambda_i ,
\]
we get
\[
 \varphi_F^* = \sum _{i\not \in I} (u_j|_{\fg^*_I} + u_j
 (\eta_0)-\lambda_j) \log (u_j|_{\fg^*_I} + u_j (\eta_0) -\lambda_j).
\]
We conclude:
\begin{theorem}
The K\"ahler form $\omega_F$
on $M_\intF$ is given by
\[
\omega_F = \sqrt{-1} \partial \bar{\partial}\left( 
(\phi_P|_{M_\intF}) ^* \sum_{j\not \in I} (
( \lambda _j - u_j (\eta_0))\log ( u_j|_{\fg^*_I} + u_j
(\eta_0)-\lambda _j)
 + u_j|_{\fg^*_I} ) \right)
\]
\end{theorem}

\subsection{Variations on the theme}

The same technique allows us to prove a variant of
(\ref{theformula}).  We keep the notation above.  Suppose that the polyhedral set $P$ is
compact.  That is, suppose that $P$ is actually a polytope.  Then
\[
\iota_\lambda (P) \subset \{ \ell \in (\R^N)^* \mid \langle \ell,
e_j\rangle \geq 0 \text{  for all }j\}
\]
is bounded.  Hence there is $R>0$ such that $\iota_\lambda (P)$ is
contained in a scaled copy $\Delta_R$ of the standard simplex:

\[
 \Delta_R = \{ \ell \in (\R^N)^* \mid \langle \ell,
e_j\rangle \geq 0 \text{  for all }j \text{ and } \sum \langle \ell,
e_j \rangle \leq R \}.
\]
Since $\Delta_1$ is the moment map image of $\C P^N$ under the
standard action of $\T^N$, it follows that $M_P$ is also a symplectic
quotient of $(\C P^N, R \omega_{FS})$ by the action of the compact
abelian Lie group $K$ defined earlier ($\omega_{FS}$ denotes the
Fubini-Study form) (c.f. \. Lemma~\ref{lem21}~(3)).  Since 
\[
 \Delta_R = \{ \ell \in (\R^N)^* \mid \langle \ell,
e_j\rangle \geq 0, \, 1\leq j \leq N,\,   \langle \ell, -\sum
e_j \rangle + R\geq 0 \},
\]
it follows from (\ref{eq-reduced-pot}) that the potential $f^*$ dual to the potential for $R \omega _{FS} $ on
$\Delta_R$ is given by
\[
f^* = \sum e_j \log e_j + (R - \sum e_j) \log (R -\sum e_j) .
\]
  Consequently the potential $f^*_\nu$ dual to the
potential on the quotient $(\C P^N/\!/_\nu K, \omega_P)$ is
\[
f^*_\nu = \sum (u_j -\lambda _j)\log (u_j -\lambda_j) + (R - \sum (u_j
-\lambda_j )) \log (R -\sum (u_j -\lambda_j)) .
\]
By Lemma~\ref{lem3.1}  the reduced 
K\"ahler form $\omega_P$ is 
\[
\omega_P = \sqrt{-1} \partial \bar{\partial} \phi^* h
\]
where
\[
h (\eta) = \langle \eta, (df_\nu)_\eta\rangle - f_\nu (\eta ).
\]
A computation similar to the ones in the previous sections gives
\begin{equation}
h = \sum \lambda_j \log (u_j - \lambda _j)
 - (R+\sum \lambda_j)\log (R -\sum (u_j -\lambda_j )).
\end{equation}

We have proved:

\begin{theorem}
Let $G$ be a torus, $P\subset \fg^*$ the polyhedral set defined by
(\ref{eqP}) which happens to be compact, $M_P = (\C P^N,
R\omega_{FS})/\!/_\nu K$ the K\"ahler $G$-space with moment map
$\phi_P:M_P \to \fg^*$ constructed in Lemma~\ref{lem21}~(3).  Then the
K\"ahler form $\omega_P$ on $\mathaccent23M_P: = \phi_P\inv
(\mathaccent23P)$ is given by:
\[
  \omega_P = \sqrt{-1} \partial \bar{\partial} \phi_P ^* \left(
\sum \lambda_j \log (u_j - \lambda _j)
 - (R+\sum \lambda_j)\log (R -\sum (u_j -\lambda_j )\right).
\]
\end{theorem}



\begin{thebibliography}{WWWW}

\bibitem[A]{Abr} M. Abreu, K\"ahler metrics on toric orbifolds, 
\emph{J.\ Diff.\ Geom.} {\bf 58} (2001), 151-187.

\bibitem[BG]{BG}  D. Burns and V. Guillemin, Potential functions and
actions of tori on K\"ahler manifolds,
\emph{Comm.\ Anal.\ Geom.} {\bf 12} (2004), no. 1-2, 281--303.

\bibitem[CDG]{CDG} D.M. Calderbank, L. David and P. Gauduchon,
 The Guillemin formula and K\"ahler metrics on toric symplectic manifolds,
\emph{J.\ Symplectic Geom.} {\bf 1} (2003), no. 4, 767--784.
 
\bibitem[D]{D} T. Delzant,  Hamiltoniens p\'eriodiques et images convexes de 
l'application moment, {\em Bull.\ Soc.\ Math.\ France} {\bf 116}
(1988), no. 3, 315--339.

\bibitem[G]{G} V. Guillemin, Kaehler structures on toric varieties,
{\em J.\ Diff.\ Geom.} {\bf 40} (1994), 285--309.

\bibitem[GS]{GS} V. Guillemin and S. Sternberg,  Geometric quantization and
multiplicities of group representations, {\em Invent.\ Math.} {\bf 67}
(1982), 515--538.

\bibitem[H]{H} P. Heinzner, e-mail communication.

\bibitem[HL]{HL} P. Heinzner and F. Loose, { Reduction of complex Hamiltonian
 $G$-spaces}, {\em Geom.\ Funct.\ Anal.} {\bf 4} (1994), 288--297.

\bibitem[HHu]{HHu} P. Heinzner and A. T. Huckleberry, K\"ahlerian
potentials and convexity properties of the moment map, {\em Invent.\
Math.} {\bf 126} (1996), 65--84.


\bibitem[HHuL]{HHuL} P. Heinzner, A. T. Huckleberry and F. Loose,
K\"ahlerian extensions of the symplectic reduction, {\em  J.\ Reine
Angew.\ Math.} {\bf  455}  (1994), 123--140.

\bibitem[MS]{MS} D. Martelli and J. Sparks, Toric geometry,
Sasaki-Einstein manifolds and a new infinite class of AdS/CFT duals,
{\tt hep-th/0411238}.

\bibitem[Sj]{Sj} R. Sjamaar, Holomorphic slices, symplectic reduction and 
multiplicities of representations, {\em Ann.\ of Math. }{\bf 141
}(1995), 87--129.


\end{thebibliography}
\end{document}